\documentclass[12pt]{amsart}
\usepackage{amssymb}

\theoremstyle{plain}
\newtheorem{theorem}{Theorem}[section]
\newtheorem{corollary}[theorem]{Corollary}
\newtheorem{lemma}[theorem]{Lemma}
\newtheorem{proposition}[theorem]{Proposition}

\theoremstyle{definition}
\newtheorem{definition}[theorem]{Definition}
\newtheorem{remark}[theorem]{Remark}

\newtheorem{example}[theorem]{Example}

\newcommand{\norm}[1]{\lVert#1\rVert}

\newcommand{\bignorm}[1]{\bigl\lVert#1\bigr\rVert}

\renewcommand{\leq}{\leqslant}
\renewcommand{\geq}{\geqslant}

\newcommand{\term}[1]{{\textit{\textbf{#1}}}}

\numberwithin{equation}{section}

\def\iA{\mathfrak A}

\def\iL{\mathcal L}
\def\iM{\mathcal M}
\def\iN{\mathcal N}
\def\iD{\mathcal D}
\def\n{{}^{(n)}}
\def\ASOT{\overline{\iA}^{\scriptscriptstyle\mathrm{SOT}}}
\def\AWOT{\overline{\iA}^{\scriptscriptstyle\mathrm{WOT}}}

\DeclareMathOperator{\Null}{Null}
\DeclareMathOperator{\Range}{Range}
\DeclareMathOperator{\Lat}{Lat}
\DeclareMathOperator{\linspan}{span}
\DeclareMathOperator{\dist}{dist}

\begin{document}

\baselineskip=18pt

\title[Strictly semi-transitive operator algebras]
      {Strictly semi-transitive\\ operator algebras}

\author[H.~P. Rosenthal, V.~G. Troitsky]{H.~P. Rosenthal and V.~G. Troitsky} 
\address{Department of Mathematics, University of Texas, Austin,
          TX 78712. USA}
\address{Department of Mathematics, University of Alberta, Edmonton,
          AB T6G\,2G1. Canada.}
\email{rosenthl@math.utexas.edu\\vtroitsky@math.ualberta.ca}

\thanks{The first author was partially supported by NSF Grant DMS 0070547}

\date{\today}

\begin{abstract}
  An algebra $\iA$ of operators on a Banach space $X$ is called
  strictly semi-transitive if for all non-zero $x,y\in X$ there exists
  an operator $A\in\iA$ such that $Ax=y$ or $Ay=x$. We show that if
  $\iA$ is norm-closed and strictly semi-transitive, then every
  $\iA$-invariant linear subspace is norm-closed. Moreover, $\Lat\iA$
  is totally and well ordered by reverse inclusion. If $X$ is complex
  and $\iA$ is transitive and strictly semi-transitive, then $\iA$ is
  WOT-dense in $\iL(X)$. It is also shown that if $\iA$ is an operator
  algebra on a complex Banach space with no invariant operator ranges,
  then~$\iA$ is WOT-dense in $\iL(X)$. This generalizes a similar
  result for Hilbert spaces proved by~Foia{\c{s}}.
\end{abstract}

\maketitle

\section{Introduction}

Throughout the paper $X$ is a real or complex Banach space, $B_X$
stands for the closed unit ball of~$X$. Let $\iA$ be a subalgebra of
$\iL(X)$, the algebra of all bounded linear operators on~$X$.  For
$x\in X$ the orbit of $x$ under $\iA$ is defined by $\iA x=\{Ax\mid
A\in\iA\}$. We say that $x_0\in X$ is \term{cyclic} for $\iA$ if $\iA
x_0$ is dense in~$X$; $x_0$ is called \term{strictly cyclic} if $\iA
x_0=X$. The algebra $\iA$ is called \term{transitive} (\term{strictly
  transitive}) if every non-zero element of $X$ is cyclic
(respectively, strictly cyclic). We write $\Lat\iA$ for the lattice of
all closed $\iA$-invariant subspaces of~$X$.  It is easy to see that
$\iA$ is transitive if and only if $\Lat\iA=\bigl\{\{0\},X\bigr\}$.
Similarly, $\iA$ is strictly transitive if and only if it has no
invariant linear subspaces other than $\{0\}$ and $X$ (by a linear
subspace we mean a subspace which is not necessarily closed). A
remarkable result due independently to Yood \cite{Yood:49} and Rickart
\cite{Rickart:50} yields that if $\iA$ is a strictly transitive
algebra of operators on a complex Banach space, then it is WOT-dense
in $\iL(X)$. We refer the reader to~\cite{Radjavi:73} for a detailed
introduction to transitive algebras.

An algebra $\iA$ is called \term{semi-transitive} if for all non-zero
$x,y\in X$ and $\varepsilon>0$ there exists an operator $A\in\iA$ such
that $\norm{Ax-y}<\varepsilon$ or $\norm{Ay-x}<\varepsilon$.  We say
that $\iA$ is \term{strictly semi-transitive} if for all non-zero
$x,y\in X$ there is an operator $A\in\iA$ with $Ax=y$ or $Ay=x$. One
can show that a unital algebra $\iA$ is semi-transitive if and only if
it is \term{unicellular}, that is, $\Lat\iA$ is totally ordered by
inclusion (see~\cite{Radjavi:73} for a treatment of unicellular
algebras).  Similarly, $\iA$ is strictly semi-transitive if and only
if all the $\iA$-invariant linear subspaces are totally ordered by
inclusion. In Section~\ref{sec:sst} we investigate the structure of
strictly semi-transitive algebras. We show that if $\iA$ is norm
closed and strictly semi-transitive, then every $\iA$-invariant linear
subspace is norm-closed and $\Lat\iA$ is well-ordered by reverse
inclusion. We deduce that a transitive strictly semi-transitive
operator algebra on a complex Banach space is WOT-dense in $\iL(X)$.
In the special case of CSL-algebras on a Hilbert space similar results
were obtained in~\cite{Hopenwasser:01,Donsig:01}.

Foia{\c{s}} proved in~\cite{Foias:72} that if $\iA$ is a WOT-closed
algebra of operators on a Hilbert space~$H$, and $\iA$ has no
invariant operator ranges, then $\iA=\iL(H)$.  In
Section~\ref{sec:op-ran} we generalize the result of Foia{\c{s}} to
complex Banach spaces using a version of Arveson's Lemma
of~\cite{Arveson:67}. For $Y\subseteq X$, we say that $Y$ is an
\term{injective operator range} if $Y=\Range\vec T$ for an injective
bounded operator $\vec T\in\iL(\iM,X)$, where $\iM$ is a closed
subspace of $X^n$ for some $n\geq 1$. We show that if $\iA$ has no
invariant injective operator ranges in $X$, then $\iA$ is WOT-dense in
$\iL(X)$.

\section{Preliminaries}

We first formulate a version of the standard lemma for the Open
Mapping Theorem, and then use it to deduce a results concerning
strictly cyclic vectors.

\begin{lemma}[Open Mapping Lemma] \label{l:open-map}
  Suppose that $Y$ is a normed space, and $T\colon X\to Y$ is a
  bounded operator. Assume that there exist $K>0$ and
  $0<\varepsilon<1$ such that
  \begin{equation}\label{eq:balls}
    B_Y\subseteq KT(B_X)+\varepsilon B_Y.
  \end{equation}
  Then $T$ is a surjective open map, and $Y$ is complete.
\end{lemma}

\begin{remark}\label{r:TB-closed}
  In particular, (\ref{eq:balls}) is satisfied if $\overline{T(B_X)}$
  has non-empty interior.
\end{remark}

\begin{proof}
  Fix $y\in B_Y$, choose $x_1\in KB_X$ such that
  $\norm{y-Tx_1}\leq\varepsilon$. Denote $y_1=y-Tx_1$, then
  $y_1\in\varepsilon B_Y$, so that there exists $x_2\in\varepsilon
  KB_X$ such that $\norm{y_1-Tx_2}\leq\varepsilon^2$, let $y_2=y_1-Tx_2$.
  Continuing, we obtain sequences $(x_n)$ in $X$ and
  $(y_n)$ in $Y$ so that for all $n$ we have
  \begin{equation}\label{eq:y_n}
    \norm{x_n}\leq\varepsilon^{n-1}K,\mbox{ and}
  \end{equation}
  \begin{equation}\label{eq:x_n}
    \norm{y_n}\leq\varepsilon^n,\mbox{ and }
      y_n=y-T\sum\limits_{j=1}^nx_j.
  \end{equation}
  Since $X$ is complete, (\ref{eq:y_n})~yields that $\sum_{j=1}^\infty
  x_j$ converges to an element $x$ of~$X$. Since $T$ is continuous,
  $Tx=y$ by~(\ref{eq:x_n}). Thus $T$ is surjective.  Moreover, it
  follows from~(\ref{eq:y_n}) that $\norm{x}\leq M$ where
  $M=K\sum_{j=0}^\infty\varepsilon^j=\frac{K}{1-\varepsilon}$, so that
  $B_Y\subseteq MT(B_X)$, hence $T$ is open. Now let $Z=X/\Null T$ and
  let $\pi\colon X\to Z$ be the canonical quotient map. Then $Z$ is
  complete and there is a one-to-one operator $\widetilde{T}\colon
  Z\to Y$ such that $\widetilde{T}\pi=T$. Since $T$ is open, the set
  $\widetilde{T}(U)=T(\pi^{-1}U)$ is open whenever $U\subseteq Z$ is
  open, so that $\widetilde{T}^{-1}$ is continuous. Hence,
  $\widetilde{T}$ is an isomorphism between $Z$ and~$Y$, so that $Y$
  is complete.
\end{proof}

\begin{proposition}\label{p:str-cyc}
  Let $\iA$ be norm closed and $x_0\in X$ with $x_0\neq 0$. The
  following are equivalent.
  \begin{enumerate}
    \item\label{i:sc} $x_0$ is strictly cyclic.
    \item\label{i:sc-unif}
      There is a constant $C>0$ so that for all $y\in B_X$ there is an
      operator $A\in\iA$ with $\norm{A}\leq C$ such that $Ax_0=y$.
    \item\label{i:sc-unif-e}
      There are $C>0$ and $0<\varepsilon<1$ so that for all $y\in B_X$
      there is an operator $A\in\iA$ with $\norm{A}\leq C$ and
      $\norm{y-Ax_0}<\varepsilon$.
  \end{enumerate}
\end{proposition}

\begin{proof}
  Define $T\colon\iA\to X$ by $T(A)=Ax_0$, then $T$ is a bounded
  linear operator from $\iA$ to~$X$. Since $\iA$ is norm closed, it
  is a Banach space. Thus,
  (\ref{i:sc})$\Rightarrow$(\ref{i:sc-unif}) follows immediately from
  the Open Mapping Theorem. The implication
  (\ref{i:sc-unif})$\Rightarrow$(\ref{i:sc-unif-e}) is
  trivial. Finally, (\ref{i:sc-unif-e})$\Rightarrow$(\ref{i:sc})
  follows from the Open Mapping Lemma.
\end{proof}

\begin{corollary}\label{c:str-cyc-open}
  If $\iA$ is norm-closed, then the set of strictly cyclic vectors for
  $\iA$ is open.
\end{corollary}

\begin{proof}
  Let $x_0$ be a strictly cyclic vector for~$\iA$, and choose $C$ as
  in Proposition~\ref{p:str-cyc}(\ref{i:sc-unif}). Let
  $0<\delta<1/C$, and suppose that $x\in X$ with
  $\norm{x-x_0}<\delta$. Now if $y\in B_X$, we choose $A\in\iA$ with
  $\norm{A}\leq C$ and $Ax_0=y$. But then
  $\norm{Ax-y}=\norm{Ax-Ax_0}\leq C\delta<1$.  Hence $x$ is strictly
  transitive by Proposition~\ref{p:str-cyc}(\ref{i:sc-unif-e}).
\end{proof}

This yields an alternate proof of a result due to
A.~Lambert.

\begin{corollary}[\cite{Lambert:71pacific}] \label{c:str-cyc-tr}
  If $X$ is a complex Banach space and $\iA$ is a transitive operator
  algebra on $X$ with a strictly cyclic vector, then $\iA$ is
  WOT-dense in $\iL(X)$.
\end{corollary}

\begin{proof}
  We may assume that $\iA$ is norm-closed by replacing $\iA$
  by~$\overline{\iA}$. Let $x_0$ be a strictly cyclic vector
  for~$\iA$.  By Corollary~\ref{c:str-cyc-open} there exists $\delta>0$
  such that $\norm{x-x_0}<\delta$ implies $x$ is strictly cyclic. Let $y\in
  X$ with $y\neq 0$. Since $\iA$ is transitive we may choose $A\in\iA$
  with $\norm{Ay-x_0}<\delta$. But then $Ay$ is strictly cyclic, and so,
  of course, $y$ is also strictly cyclic. Thus, $\iA$ is strictly
  transitive and so is WOT-dense in $\iL(X)$.
\end{proof}

The next result shows that transitive algebras always
have operators which are almost zero on prescribed vectors.

\begin{proposition}\label{p:almost-zero}
  Let $\iA$ be a transitive operator algebra on a complex Banach space
  $X$, then for all $x\in X$ and $\varepsilon>0$ there is an $A\in\iA$
  with $\norm{A}=1$ and $\norm{Ax}<\varepsilon$.
\end{proposition}

\begin{proof}
  Without loss of generality, $\iA$ is WOT-closed. Indeed, suppose
  that the conclusion holds for $\AWOT\!\!$, show that it also holds
  for $\iA$. Fix $x\in X$ and $0<\varepsilon<1$, and let
  $\delta=\varepsilon/4$. There exists $A\in\AWOT$ with $\norm{A}=1$
  and $\norm{Ax}<\delta$. Then $\norm{Ay}\geq 1-\delta$ for some $y\in
  X$ with $\norm{y}=1$. Since $\AWOT\!\!=\ASOT\!\!$, there exists
  $B\in\iA$ such that $\bignorm{(A-B)x}<\delta$ and
  $\bignorm{(A-B)y}<\delta$. Then $\norm{Bx}<2\delta$ and
  $\norm{By}\geq 1-2\delta$, so that $\norm{B}\geq 1-2\delta$. Let
  $C=\frac{B}{\norm{B}}$, then
  $\norm{Cx}<\frac{2\delta}{1-2\delta}<\varepsilon$.

  If $\iA=\iL(X)$ then the conclusion is trivially satisfied. Suppose
  that $\iA$ is a WOT-closed proper subalgebra of $\iL(X)$. Were the
  conclusion false, we could choose $x$ of norm one and
  $\varepsilon>0$ so that $\norm{Ax}>\varepsilon$ whenever $A\in\iA$
  and $\norm{A}=1$. But then the operator $T\colon\iA\to X$ defined by
  $T(A)=Ax$ for $A\in\iA$ is an isomorphism. In particular, $\iA
  x$ is norm-closed being isomorphic to~$\iA$. Then $\iA x=X$
  since $\iA$ is transitive. This contradicts
  Corollary~\ref{c:str-cyc-tr}.
\end{proof}

The following generalization of this result is due to Jiaosheng
Jiang, and we are grateful to him for giving us permission to present
his proof here.

\begin{theorem}[J.~Jiang] \label{t:jiang}
  Let $\iA$ be a commutative transitive operator algebra on a complex
  Banach space. Then for all $x_0,x_1,\dots x_n\in X$ and
  $\varepsilon>0$ there exists an operator $A\in\iA$ with $\norm{A}=1$
  and $\norm{Ax_i}<\varepsilon$ for all $0\leq i\leq n$.
\end{theorem}

\begin{proof}
  We may assume without loss of generality that $x_0\neq 0$. 
  By the transitivity of $\iA$, for each $1\leq i\leq n$ we may choose
  $B_i\in\iA$ so that $\norm{B_ix_0-x_i}\leq\varepsilon/2$.
  Let $\delta=\min\bigl\{\varepsilon,
    \frac{\varepsilon}{2\cdot\max_i\norm{B_i}}\bigr\}$.
  By Proposition~\ref{p:almost-zero}, we may choose $A\in\iA$ with
  $\norm{A}=1$ and $\norm{Ax_0}<\delta$. Then
  $\norm{Ax_0}<\varepsilon$ and for all $1\leq i\leq n$,
  \begin{multline*}
    \norm{Ax_i}\leq\bignorm{A(x_i-B_ix_0)}+\norm{AB_ix_0}\\
    \leq\norm{x_i-B_ix_0}+\norm{B_iAx_0}
    \leq\varepsilon/2+\norm{B_i}\cdot\delta<\varepsilon.
  \end{multline*}
\end{proof}

We do not know if the conclusion of this proposition holds for general
transitive algebras. The following result gives a partial answer.
Recall that $\iA$ is said to be \term{$n$-transitive} for some $n\geq
1$ if every linearly independent $n$-tuple in $X^n$ is cyclic for
the algebra $\iA\n=\bigl\{A\oplus\dots\oplus A\mid
A\in\iA\bigr\}$.

\begin{proposition}
  Let $\iA$ be an operator algebra on a complex Banach space. If $\iA$
  is $n$-transitive for some $n\geq 1$, then for all $x_1,\dots,x_n\in
  X$ and $\varepsilon>0$ there exists an operator $A\in\iA$ with
  $\norm{A}=1$ and $\norm{Ax_i}<\varepsilon$ for all~$i$.
\end{proposition}

\begin{proof}
  As in the proof of Proposition~\ref{p:almost-zero}, we may assume that
  $\iA$ is a WOT-closed proper subalgebra of $\iL(X)$. Were the
  conclusion false, we could find $x_1,\dots,x_n\in X$ and
  $\varepsilon>0$ such that $\max\limits_{1\leq i\leq
    n}\norm{Ax_i}\geq\varepsilon$ for all $A\in\iA$ with $\norm{A}=1$.
  Without loss of generality, we may assume that $x_1,\dots,x_k$ are
  linearly independent for some $k\leq n$, and
  $x_j\in\linspan\{x_1,\dots,x_k\}$ whenever $k<j\leq n$. It follows
  that there is a $C>1$ such that
  \begin{displaymath}
    \max\limits_{1\leq i\leq n}\norm{Ax_i}\leq
    C\cdot\max\limits_{1\leq i\leq k}\norm{Ax_i}
  \end{displaymath}
  whenever $A\in\iA$ with $\norm{A}=1$, so that $\max\limits_{1\leq
    i\leq k}\norm{Ax_i}\geq\varepsilon/C$. Then the map $T\colon\iA\to
  X^k$ defined by $T(A)=(Ax_1,\dots,Ax_k)$ for $A\in\iA$ is an
  isomorphism, and hence $T(\iA)$ is closed. Since $\iA$ is
  $n$-transitive, it is also $k$-transitive. Then $T(\iA)=X^k$ because
  $x_1,\dots,x_k$ are linearly independent. This implies, in
  particular, that $x_1$ is strictly cyclic for $\iA$, and, since
  $\iA$ is transitive, $\iA$ is WOT-dense in $\iL(X)$ by
  Corollary~\ref{c:str-cyc-tr}, a contradiction.
\end{proof}

Recall that $T\in\iL(X,Y)$ is called a \term{semi-embedding} if 
it is one-to-one and $T(B_X)$ is closed. The following result is
certainly known, and we include a proof for the sake of completeness.

\begin{proposition} \label{p:semi-emb}
  Let $Y$ be a normed space, and $T\colon Y\to X$ a
  semi-embedding. Then $Y$ is complete.
\end{proposition}

\begin{proof}
  Let $S\colon\Range(T)\to Y$ be the inverse of~$T$, and let
  $\widetilde{Y}$ be the completion of~$Y$.  Then $T$ may be extended
  to a continuous operator $\widetilde{T}\colon\widetilde{Y}\to X$.
  Since $T(B_Y)$ is closed in $X$ we have
  $\widetilde{T}(B_{\widetilde{Y}})=T(B_Y)$, and hence
  $\Range\widetilde{T}=\Range T$. Consider the operator
  $S\widetilde{T}\colon\widetilde{Y}\to Y$. If $y\in Y$ then
  $S\widetilde{T}y=STy=y$.  Moreover,
  $S\widetilde{T}(B_{\widetilde{Y}})=S\bigl(T(B_Y)\bigr)= B_Y\subseteq
  B_{\widetilde{Y}}$.  It follows that $S\widetilde{T}$ is a bounded
  projection of $\widetilde{Y}$ onto~$Y$, and thus $Y$ is complemented
  in~$\widetilde{Y}$, so that it is closed in~$\widetilde{Y}$, hence
  $\widetilde{Y}=Y$.
\end{proof}

We are now prepared for the final result of this section.  For
$W\subseteq X$ and $\varepsilon>0$, let $W_\varepsilon=\bigl\{x\in
X\mid\dist(x,W)\leq\varepsilon\bigr\}$. The continuity of the map
$x\mapsto\dist(x,W)$ from $X$ to $\mathbb R$ implies that
$W_\varepsilon$ is a closed set.

\begin{lemma} \label{l:no-balls}
  Let $X$ be a Banach space and $W$ a bounded circled closed convex
  subset of $X$ with void interior. Then for any $\varepsilon>0$ the set
  $W_\varepsilon$ contains no ball of radius larger than
  $\varepsilon$.
\end{lemma}

\begin{proof}
  Let $Y=\linspan W=\bigcup_{n=1}^\infty nW$. Then $Y$ equipped with
  the norm given by the Minkowski functional of $W$ is a normed space.
  Let $T\colon Y\to X$ be the inclusion operator. Then $T$ is bounded,
  one-to-one, and $T(B_Y)=W$. Hence, $T$ is a semi-embedding, so by
  Proposition~\ref{p:semi-emb} $Y$ is complete.
  
  Clearly, $W_\varepsilon$ is also circled and convex.  Suppose that
  $W_\varepsilon$ contains a ball of radius $r>\varepsilon$ centered
  at~$x_0$. Since $W_\varepsilon$ is circled, it also contains the
  ball of radius $r$ centered at $-x_0$, and since $W_\varepsilon$ is
  convex, it contains the ball of radius $r$ centered at the origin.
  Indeed, if $\norm{z}\leq r$ then $x_0+z$ and $-x_0+z$ belong
  to~$W_\varepsilon$, so that $z=\frac{(x_0+z)+(-x_0+z)}{2}\in
  W_\varepsilon$. Now let $\varepsilon<\varepsilon'<r$. Then
  $rB_X\subseteq W_\varepsilon\subseteq W+\varepsilon'B_X$, thus,
  $B_X\subseteq\frac{1}{r}W+\frac{\varepsilon'}{r}B_X$. Since
  $\frac{\varepsilon'}{r}<1$, we conclude by by Lemma~\ref{l:open-map}
  that $T$ is surjective. This yields $X=Y=\bigcup_{n=1}^\infty nW$,
  which contradicts the Baire Category Theorem since $W$ has void
  interior.
\end{proof}

\section{The structure of strictly semi-transitive algebras}
\label{sec:sst}

A subalgebra $\iA$ of $\iL(X)$ is called \term{strictly
semi-transitive} if for all non-zero $x,y\in X$ there exists an
operator $A\in\iA$ so that $Ax=y$ or $Ay=x$.
\begin{proposition}\label{p:sst-equiv}
  Let $\iA$ be a unital subalgebra of $\iL(X)$. Then the following are
  equivalent.
  \begin{enumerate}
    \item\label{i:sst} $\iA$ is strictly semi-transitive;
    \item\label{i:subsp-ord} All the $\iA$-invariant linear
      subspaces of $X$ are totally ordered by inclusion;
    \item\label{i:orb-ord} All the orbits are totally ordered by
      inclusion.
  \end{enumerate}
\end{proposition}

\begin{proof}
  (\ref{i:sst})$\Rightarrow$(\ref{i:subsp-ord}). Let $Y$ and $Z$ be
  two $\iA$-invariant linear subspaces of~$X$, and suppose that
  $Z\nsubseteq Y$. Show that $Y\subseteq Z$. Choose $z\in Z\setminus
  Y$ and $y\in Y$. There is no $T\in\iA$ with $Ty=z$, for otherwise
  $z\in Y$ since $Y$ is $\iA$-invariant. Hence there must be a
  $T\in\iA$ with $Tz=y$, and so $y\in Z$ since $Z$ is $\iA$-invariant.

  (\ref{i:subsp-ord})$\Rightarrow$(\ref{i:orb-ord}) is trivial.

  (\ref{i:orb-ord})$\Rightarrow$(\ref{i:sst}). Let $x,y\in X$. Then
  $\iA x\subseteq\iA y$ or $\iA y\subseteq\iA x$. Since $\iA$ is
  unital, it follows that $x\in\iA y$ or $y\in\iA x$, hence $\iA$ is
  strictly semi-transitive.
\end{proof}

\begin{remark}\label{r:sst-no-unit}
  The implications
  (\ref{i:sst})$\Rightarrow$(\ref{i:subsp-ord})$\Rightarrow$(\ref{i:orb-ord})
  are still valid for non-unital algebras.
\end{remark}

The next result is fundamental for our development. Recall that a set
is \term{residual} if its complement is of first category; residual
subsets of $X$ are dense by Baire Category Theorem.

\begin{lemma} \label{l:residual}
  Let $(W_n)$ be a sequence of circled bounded closed convex sets,
  each with void interior, such that $W_n\subseteq W_{n+1}$ and
  $\bigcup\limits_{n=1}^\infty W_n$ is dense in $X$, and let
  $(\varepsilon_n)$ and $(\delta_n)$ be sequences of positive reals
  tending to zero. Let
  \begin{eqnarray*}
    G_1&=&\bigl\{x\in X\mid\text{ for all }\lambda>0,\,
          \dist(\lambda x,W_n)>\varepsilon_n
          \text{ for infinitely many }n\},\\
    G_2&=&\bigl\{x\in X\mid\text{ for all }\lambda>0,\,
          \dist(\lambda x,W_n)<\delta_n
          \text{ for infinitely many }n\}.
  \end{eqnarray*}
  Then $G_1\cap G_2$ is residual.
\end{lemma}

\begin{proof}
  Note that $x\in G_1$ if and only if
  \begin{equation}\label{eq:G1-real}
    \forall\,\lambda>0\quad\forall\,m\in\mathbb N\quad\exists\, n\geq m
    \text{ such that }\dist(\lambda  x,W_n)>\varepsilon_n.
  \end{equation}
  In particular,
  \begin{equation}\label{eq:G1-natural}
    \forall\,k\in\mathbb N\quad\forall\,m\in\mathbb N\quad\exists\,n\geq m
    \text{ such that }\dist\bigl(\tfrac{1}{k}x,W_n\bigr)>\varepsilon_n.
  \end{equation}
  We claim that (\ref{eq:G1-real}) and~(\ref{eq:G1-natural}) are, in
  fact, equivalent. Indeed, suppose that (\ref{eq:G1-natural}) holds.
  Given $\lambda>0$ and $m\in\mathbb N$, find a positive integer $k$ such
  that $\frac{1}{k}<\lambda$, then there exists $n\geq m$ such that
  $\dist(\frac{1}{k}x,W_n)>\varepsilon_n$. Since $W_n$ is convex and circled,
  it follows that $\dist(\lambda
  x,W_n)\geq\dist(\frac{1}{k}x,W_n)>\varepsilon_n$. 
  
  It follows from~(\ref{eq:G1-natural}) that $x\notin G_1$ if and only
  if there exist $k,m\in\mathbb N$ such that
  $\frac{1}{k}x\in(W_n)_{\varepsilon_n}$ for all $n\geq m$. Therefore,
  \begin{displaymath}
    \sim\! G_1
        =\bigcup\limits_{k,m=1}^\infty
          \bigcap\limits_{n=m}^\infty k(W_n)_{\varepsilon_n}
        =\bigcup\limits_{k,m=1}^\infty
        k\bigcap\limits_{n=m}^\infty(W_n)_{\varepsilon_n}.
  \end{displaymath}
  It follows from Lemma~\ref{l:no-balls} that $(W_n)_{\varepsilon_n}$
  contains no balls of radius larger than $\varepsilon_n$. Therefore,
  for every $m$ the set
  $\bigcap_{n=m}^\infty(W_n)_{\varepsilon_n}$ has void interior, so that
  $\sim\! G_1$ is of first category.
  
  Similarly, $x\in G_2$ if and only if
  \begin{equation}\label{eq:G2-natural}
    \forall\,k\in\mathbb N\quad\forall\,m\in\mathbb N\quad\exists\,n\geq m
    \text{ such that }\dist(kx,W_n)<\delta_n.
  \end{equation}
  Indeed, if $x\in G_2$ then (\ref{eq:G2-natural}) is satisfied
  trivially. Conversely, suppose that $x$
  satisfies~(\ref{eq:G2-natural}), let $\lambda>0$ and $m\in\mathbb
  N$. Let $k\in\mathbb N$ be such that $\lambda\leq k$.
  By~(\ref{eq:G2-natural}) there exists $n\geq m$ such that
  $\dist(kx,W_n)<\delta_n$. But since $W_n$ is convex and circled we
  have $\dist(\lambda x,W_n)\leq\dist(kx,W_n)<\delta_n$. Note that
  $\dist(kx,W_n)<\delta_n$ is equivalent to
  $x\in\frac{1}{k}(W_n+\delta_nB_X^\circ)$, where $B_X^\circ$ stands
  for the open unit ball of~$X$. Thus,
  \begin{displaymath}
    G_2=\bigcap\limits_{k,m=1}^\infty
        \tfrac{1}{k}\bigcup\limits_{n=m}^\infty
        (W_n+\delta_nB_X^\circ).
   \end{displaymath}
   Since $\bigcup_{n=m}^\infty(W_n+\delta_nB_X^\circ)$ is open and
   dense in $X$ for every~$m$, then by Baire Theorem $G_2$ is
   residual.  Since $G_1$ and $G_2$ are both residual, so is their
   intersection.
\end{proof}

The next result is our main lemma. It shows that if $W_n$'s are as in
the previous lemma, then there are always $u$ and $v$ in $X$ which
``see'' the $W_n$'s very differently.

\begin{lemma}\label{l:see-diff}
  Let $(W_n)$ be a nested increasing sequence of closed convex bounded
  circled sets, each with void interior, so that
  $Y=\bigcup_{n=1}^\infty W_n$ is dense in~$X$. Then given any $u\in
  X\setminus Y$, there is a vector $v\in X\setminus Y$ so that
  \begin{equation}\label{eq:limsups}
     \limsup\limits_{n\to\infty}\frac{\dist(v,W_n)}{\dist(u,W_n)}=\infty
     \quad\mbox{ and }\quad
     \liminf\limits_{n\to\infty}\frac{\dist(v,W_n)}{\dist(u,W_n)}=0.
  \end{equation}
\end{lemma}

\begin{proof}
  Note that $Y$ is a linear subspace of $X$ of first category, so
  $Y\neq X$ by the Baire Category Theorem. Let $\tau_n=\dist(u,W_n)$,
  then the sequence $(\tau_n)$ is decreasing and
  $\lim\limits_{n\to\infty}\tau_n=0$.  Now let $G_1$ and $G_2$ be as
  in Lemma~\ref{l:residual} with $\varepsilon_n=\sqrt{\tau_n}$ and
  $\delta_n=\tau_n/n$, then $G_1\cap G_2$ is residual, hence
  non-empty. Let $v\in G_1\cap G_2$, then
  \begin{displaymath}
    \frac{\dist(v,W_n)}{\dist(u,W_n)}>
    \frac{\varepsilon_n}{\tau_n}=
    \frac{1}{\sqrt{\tau_n}}
  \end{displaymath}
  for infinitely many~$n$, proving the first equality
  in~(\ref{eq:limsups}), and
  \begin{displaymath}
    \frac{\dist(v,W_n)}{\dist(u,W_n)}<
    \frac{\delta_n}{\tau_n}=\frac{1}{n}
  \end{displaymath}
  for infinitely many~$n$, proving the second equality.
\end{proof}

\begin{remark}\label{r:resid-homog}
  Notice that every $v$ in $G_1\cap G_2$, which is a residual set
  in~$X$, satisfies~(\ref{eq:limsups}). Furthermore, by definition the
  sets $G_1$ and $G_2$ are positively homogeneous, that is, with every
  vector they contain the entire ray passing through that vector, so
  that $G_1\cap G_2$ is also positively homogeneous.
\end{remark}

We are now ready for our first main structural result for strictly
semi-transitive algebras.

\begin{theorem}\label{t:closed-orb}
  Let $\iA$ be norm-closed and strictly semi-transitive. Then $\iA x$
  is closed for every $x\in X$.
\end{theorem}

\begin{proof}
  Let $X'=\overline{\iA x}$.  Define $T\colon\iA\to X'$ by $T(A)=Ax$
  for $A\in\iA$, and let
  \begin{displaymath}
    W_n=n\overline{T(B_{\iA})}=
        \overline{\bigl\{Ax\mid A\in\iA,\,\norm{A}\leq n\bigr\}}
  \end{displaymath}
  for all~$n$. Of course, $W_n=nW_1$ and $W_n\subset X'$ for every~$n$. It
  suffices to prove that $W_1$ has non-void relative interior in $X'$
  because in this case Remark~\ref{r:TB-closed} would imply that
  $X'=T(\iA)=\iA x$. Suppose not. Then, of course, each $W_n$ is a
  closed convex bounded circled set, with void relative interior
  in~$X'$. Pick any $u$ in $X'\setminus\bigcup_{n=1}^\infty W_n$. It follows from
  Lemma~\ref{l:see-diff} that we may choose $v$ in $X'\setminus\bigcup_{n=1}^\infty
  W_n$ satisfying~(\ref{eq:limsups}).  Since $\iA$ is strictly
  semi-transitive, there is an operator $A\in\iA$ so that $Au=v$ or
  $Av=u$.  It follows from Remark~\ref{r:resid-homog} that by scaling
  $v$ and $A$ we can also assume without loss of generality that
  $\norm{A}=1$.  Then $A(W_n)\subseteq W_n$ for every~$n$. If $Au=v$, then
  \begin{displaymath}
    \dist(v,W_n)\leq
    \dist(Au,AW_n)\leq
    \dist(u,W_n),
  \end{displaymath}
  for every $n$, but this contradicts the first equality
  in~(\ref{eq:limsups}). Similarly, the second equality
  in~(\ref{eq:limsups}) is violated if $Av=u$. This contradiction
  completes the proof.
\end{proof}

\begin{theorem}\label{t:sst-trans-dense}
  Let $\iA$ be a strictly semi-transitive operator algebra on a
  complex Banach space~$X$. If $\iA$ is transitive, then it is
  WOT-dense in $\iL(X)$.
\end{theorem}

\begin{proof}
  We may assume without loss of generality that $\iA$ is norm-closed.
  But then $\iA$ is strictly transitive by Theorem~\ref{t:closed-orb},
  hence WOT-dense in $\iL(X)$.
\end{proof}

We now arrive at a second basic structural result.

\begin{theorem} \label{t:lat-well-ord}
  Let $\iA$ is strictly semi-transitive. Then $\Lat\iA$ is
  well-ordered by reverse inclusion. That is, every non-empty subset
  of $\Lat\iA$ has a maximal element.
\end{theorem}

\begin{proof}
  We may without loss of generality assume that $\iA$ is norm closed.
  Suppose that $\mathcal W$ is a non-empty subset of $\Lat\iA$ with no
  maximal element. Then we can find an infinite sequence
  $X_1\subsetneqq X_2\subsetneqq X_3\subsetneqq\ldots$ in~$\mathcal W$.
  Let $Y=\bigcup_{n=1}^\infty X_n$. Then $Y$ is not closed; indeed each $X_n$
  is nowhere dense in the induced topology of~$\overline{Y}$, so that
  $Y$ is of first category in~$\overline{Y}$. Furthermore,
  $\overline{Y}$ is $\iA$-invariant (i.e., $\overline{Y}\in\Lat\iA$),
  and $\iA$ reduced to $\overline{Y}$ is again a strictly
  semi-transitive algebra.  Set $W_n=nB_{X_n}$ for all~$n$. Again by
  Lemma~\ref{l:see-diff} we may choose $u,v\in\overline{Y}\setminus Y$
  satisfying~(\ref{eq:limsups}). Again, choose $A\in\iA$ of norm one so
  that $Au=v$ or $Av=u$. Then $A(W_n)\subseteq W_n$ because $X_n$ is an
  invariant subspace for~$\iA$. Now the rest of the argument for
  Theorem~\ref{t:closed-orb} yields a contradiction.
\end{proof}

The next result yields the surprising fact that every $\iA$-invariant
linear subspace is in $\Lat\iA$.

\begin{corollary}\label{c:inv-closed}
  Let $\iA$ be unital, norm-closed and strictly semi-transitive. Then
  every $\iA$-invariant linear subspace $Y$ of $X$ is closed.
  Furthermore, $Y=\iA x$ for some $x\in X$.
\end{corollary}

\begin{proof}
  Let $\mathcal W$ be the set of all orbits contained in~$Y$, i.e.,
  $\mathcal W=\{\iA x\mid\iA x\subseteq Y\}$. Note that $\mathcal
  W\subseteq\Lat\iA$ by Theorem~\ref{t:closed-orb}. Then
  Theorem~\ref{t:lat-well-ord} yields that $\mathcal W$ has a maximal
  element, say $\iA x$. We claim that $Y=\iA x$. Suppose not, then
  there exists $y\in Y\setminus\iA x$. Then $\iA x\subseteq\iA
  y\subseteq Y$ by Proposition~\ref{p:sst-equiv}, and $\iA x\neq\iA y$
  because $y\in\iA y\setminus\iA x$. But this contradicts the
  maximality of $\iA x$.  Hence $Y=\iA x$, so that $Y$ is closed by
  Theorem~\ref{t:closed-orb}.
\end{proof}

The following result may alternatively be deduced
from~\cite[Theorem~4.4]{Radjavi:73} and Theorem~\ref{t:closed-orb}. We
present a simple direct proof of it.

\begin{corollary}
  Let $\iA$ be norm-closed and strictly semi-transitive. Then $\iA$
  has (a residual set of) strictly cyclic vectors.
\end{corollary}

\begin{proof}
  If $\iA$ is transitive, then every non-zero vector is strictly
  cyclic by Theorem~\ref{t:closed-orb}. Otherwise, $\iA$ has a
  maximal proper closed invariant subspace $Y$ by
  Theorem~\ref{t:lat-well-ord}. Then $X\setminus Y$ is residual. But
  if $x\in X\setminus Y$, then $\iA x=X$ by Theorem~\ref{t:closed-orb}
  and the maximality of~$Y$.
\end{proof}

We conclude this section with an explicit description of $\Lat\iA$
when $\iA$ a strictly semi-transitive algebra acting on a separable
space, and a description of all ``full'' strictly semi-transitive
algebras.

\begin{corollary}\label{c:count-ordinal}
  Let $X$ be separable and $\iA$ be strictly semi-transitive. Then
  there exists a countable ordinal $\eta$ and a family of closed
  subspaces $(X_\alpha)_{\alpha\leq\eta}$ such that
  \begin{enumerate}
  \item\label{i:nested} $X_0=X$, $X_{\eta}=\{0\}$, and
    $X_{\alpha+1}\subsetneqq X_{\alpha}$ for all $\alpha<\eta$;
  \item $\Lat\iA=\{X_\alpha\mid\alpha\leq\eta\}$. 
  \end{enumerate}
\end{corollary}

\begin{proof}
  This holds for arbitrary $X$ by Theorem~\ref{t:lat-well-ord}, except
  that $\eta$ need not be countable. However, if $X$ is assumed to be
  separable, it follows that $\eta$ must be countable, otherwise $X$
  would have an uncountable strictly decreasing chain of closed
  subspaces, which is impossible.
\end{proof}

\begin{remark}
  Assume that $\eta$ and $X_\alpha$'s are as in the above
  statement. It follows that
  \begin{gather}\label{eq:sequence}
    \begin{split}
      \textit{if $\alpha_1<\alpha_2<\dots$ with $\alpha_n<\eta$ for all $n$}\\
      \textit{and $\alpha_n\to\alpha$, then $X_\alpha=\bigcap_{n=1}^\infty
      X_{\alpha_n}$.}
    \end{split}
  \end{gather}
  Indeed, let
  $Y=\bigcap_{n=1}^\infty X_{\alpha_n}$, then $Y\in\Lat\iA$, hence
  $Y=X_\beta$ for some~$\beta$. But then $\beta\geq\alpha_n$ for all
  $n$, hence $\beta\geq\alpha$. This yields $X_\beta\subseteq
  X_\alpha$, but $X_\alpha\subseteq X_{\alpha_n}$ for all $n$, so that
  \begin{displaymath}
    \bigcap\limits_{n=1}^\infty X_{\alpha_n}=X_\beta\subseteq X_\alpha
    \subseteq\bigcap\limits_{n=1}^\infty X_{\alpha_n},
  \end{displaymath}
  hence $X_\beta=X_\alpha$.
\end{remark}

The above remark motivates our final result of this section, which is
a partial converse to Theorem~\ref{t:lat-well-ord}.

\begin{proposition}
  Let $X$ be a separable Banach space, $\eta$ a countable ordinal with
  $\eta>1$, and $(X_\alpha)_{\alpha\leq\eta}$ a family of closed
  subspaces of $X$ satisfying~(\ref{i:nested}) of
  Corollary~\ref{c:count-ordinal} and~(\ref{eq:sequence}). Let
  \begin{displaymath}
    \iA=\{T\in\iL(X)\mid TX_\alpha\subseteq X_\alpha\text{ for all
    }\alpha\leq\eta\}.
  \end{displaymath}
  Then $\iA$ is strictly semi-transitive and
  $\Lat\iA=\{X_\alpha\mid\alpha\leq\eta\}$.
\end{proposition}

\begin{proof}
  Let $x,y\in X\setminus\{0\}$. Put
  \begin{displaymath}
    \alpha=\sup\{\gamma\leq\eta\mid x\in X_\gamma\}\qquad\mbox{ and }\qquad
    \beta=\sup\{\gamma\leq\eta\mid y\in X_\gamma\}
  \end{displaymath}
  It follows from~(\ref{eq:sequence}) that $x\in X_\alpha$ and $y\in
  X_\beta$. Without loss of generality, $\alpha\leq\beta$ and, of
  course, $\beta<\eta$. Then by the Hahn Banach Theorem, there is an
  $f\in X^*$ with $X_{\beta+1}\subseteq\ker f$ and $f(x)=1$. Define
  $T\in\iL(X)$ by $Tu=f(u)y$ for all $u\in X$. Evidently $Tx=y$. If
  $\gamma\geq\beta+1$ then $X_\gamma\subseteq X_{\beta+1}\subseteq\ker
  f$, so $TX_\gamma=\{0\}$, and if $\gamma\leq\beta$, then
  $X_\beta\subseteq X_\gamma$ and, of course, $TX=[y]\subseteq
  X_\gamma$. In either case, $TX_\gamma\subseteq X_\gamma$. Thus,
  $T\in\iA$, and we have shown that $\iA$ is strictly semi-transitive.
  
  It follows by definition that $\{X_\alpha\mid\alpha\leq\eta\}\subseteq\Lat\iA$, so we
  must prove that, conversely, if $Y\in\Lat\iA$, then $Y=X_\alpha$ for
  some~$\alpha$. Assume $Y\neq\{0\}$ or $X$, and let $\alpha$ be the greatest
  such that $Y\subseteq X_\alpha$. We claim that $Y=X_\alpha$. Suppose this were
  false. Now since $Y\subsetneqq X_\alpha$, it follows from
  Proposition~\ref{p:sst-equiv} and the definition of $\alpha$ that
  $X_{\alpha+1}\subsetneqq Y$. Now we may choose $z\in X_\alpha\setminus Y$ and $y\in
  Y\setminus X_{\alpha+1}$. But then by our initial discussion, since neither $z$
  nor $y$ are in $X_{\alpha+1}$, yet both belong to $X_\alpha$, there is a
  $T\in\iA$ with $Ty=z$. But then $TY\nsubseteq Y$, contradicting he
  assumption that $Y\in\Lat\iA$.
\end{proof}

\begin{example}
  Let $X=\ell_p$ for $1\leq p<+\infty$ and let $(e_i)$ be the standard
  basis of $\ell_p$. Let $\iA$ be the set of all the bounded operators
  $A$ on $\ell_p$ such that
  \begin{enumerate}
    \item the matrix of $A$ is lower triangular, that is,
      $Ae_n\in[e_i]_{i=n}^\infty$ for every $n$, and
    \item the matrix of $A$ contains only finitely many non-zero
      columns or, equivalently, there exists a positive integer $n$
      such that $Ae_i=0$ for all $i\geq n$.
  \end{enumerate}
  It is easy to see that $\iA$ is an algebra.  Show that $\iA$ is
  strictly semi-transitive. Given $x=(x_i)$ and $y=(y_i)$ in
  $\ell_p$, let $x_k$ be the first non-zero component of $x$ and
  $y_m$ the first non-zero component of $y$. Without loss of
  generality, $k\leq m$. Define an operator $A$ as follows:
  $Ae_k=\frac{1}{x_k}y$ and $Ae_i=0$ for all $i\neq k$. Then
  $A\in\iA$ and $Ax=y$.
  
  Notice that this algebra is not norm closed. However, it follows
  that the algebra of all lower-triangular compact operators and the
  algebra of all lower-triangular bounded operators are strictly
  semi-transitive.
  
  On the other hand, consider the algebra of all upper-triangular
  bounded operators. The linear subspace of all sequences with
  finitely many non-zero entries is invariant under the algebra, but
  it is not closed in~$\ell_p$. It follows from
  Corollary~\ref{c:inv-closed} that the upper-triangular algebra is
  not strictly semi-transitive. However, it is easy to see that this
  algebra is unicellular, hence semi-transitive.
\end{example}

\section{Operator ranges in Banach spaces}
\label{sec:op-ran}

A subspace $Y$ of a Hilbert space $H$ is termed an operator range by
Foia{\c{s}}~\cite{Foias:72} if there exists an operator $T\in\iL(H)$
with $Y=T(H)$. We generalize this to arbitrary Banach spaces $X$ as
follows.

\begin{definition}\label{d:op-ran}
  Let $X$ be a Banach space and $Y$ a linear subspace of~$X$.  Then
  $Y$ is called an \term{operator range} if there is a closed subspace
  $\iM\subseteq X^n$ for some $n\geq 1$ and a bounded linear operator $\vec
  T\colon\iM\to X$ with $Y=\Range\vec T$. In this case, $Y$ is called an
  \term{operator range of order~$n$}. Finally, $Y$ is called an
  \term{injective operator range} if $\vec T$ may be chosen
  one-to-one.
\end{definition}

Of course, every closed subspace of $X$ is an operator range of order
one. It is proved in~\cite{Foias:72} (see
also~\cite[Theorem~8.9]{Radjavi:73}) that if $\iA$ is a WOT-closed
algebra of operators on a Hilbert space with no non-trivial invariant
operator ranges, then $\iA=\iL(H)$. The main result of this section
generalizes this fact.

\begin{theorem}\label{t:inj-op-ran}
  Let $\iA$ be a subalgebra of $\iL(X)$ for a complex Banach space
  $X$, and $n\geq 1$. Then either $\iA$ has a non-trivial invariant
  injective operator range of order $n$, or $\iA$ is $n$-transitive.
\end{theorem}

Note that $\iA$ is WOT-dense in $\iL(X)$ if and only if $\iA$ is
$n$-transitive for all~$n$, c.f.~\cite[Theorem~7.1]{Radjavi:73}.
This yields the following consequence.

\begin{corollary}\label{c:WOT-dense}
  Let $\iA$ be a subalgebra of $\iL(X)$ for a complex Banach space $X$
  such that $\iA$ has no non-trivial invariant injective operator
  ranges.  Then $\iA$ is WOT-dense in $\iL(X)$.
\end{corollary}

\begin{remark}
  It is easily seen that every operator range in a Hilbert space is
  injective of order one. It seems we need the concept of operator
  range of order~$n$, however, to obtain a general result. Also, we
  don't know if every operator range in a Banach space is also an
  injective operator range. This question was posed by P.~Rosenthal.
\end{remark}

Before proving the theorem, we need some preliminary results, which
also hold in real Banach spaces. We start with the following simple
consequences of the definition of an operator range, given
in~\cite{Foias:72} in a Hilbert space case.

\begin{proposition}\label{p:op-ran-lat}
  Let $Y_1$ and $Y_2$ be operator ranges in~$X$. Then $Y_1+Y_2$ and
  $Y_1\cap Y_2$ are operator ranges.
\end{proposition}

\begin{proof}
  For $i=1,2$ choose $n_i$, $\iM_i$ a closed subspace of $X^{n_i}$,
  and $T_i\colon\iM_i\to X$ bounded linear with $Y_i=\Range
  T_i$. Now if
  \begin{displaymath}
    \iM=\iM_1\oplus\iM_2\subseteq X^{n_1}\oplus X^{n_2}=X^{n_1+n_2},   
  \end{displaymath}
  and $T=T_1\oplus T_2$, then $T(\iM)=Y_1+Y_2$, hence $Y_1+Y_2$ is an
  operator range. The argument for $Y_1\cap Y_2$ is not quite obvious.
  Let $$\mathcal W=\{(w_1,w_2)\in\iM\mid T_1w_1=T_2w_2\}.$$
  Evidently, $\mathcal W$ is a linear subspace of~$\iM$.
  If $(u_n,v_n)\to(u,v)$
  with $(u_n,v_n)\in\mathcal W$ for all~$n$, then
  $T_1u=\lim\limits_{n\to\infty}T_1u_n=\lim\limits_{n\to\infty}T_2v_n=T_2v$, so
  that $(u,v)\in\mathcal W$, hence $\mathcal W$ is closed. Show that
  $T(\mathcal W)=Y_1\cap Y_2$.  Indeed, if $(w_1,w_2)\in\mathcal W$ then
  $T(w_1,w_2)=T_1w_1+T_2w_2=2T_1w_1=2T_2w_2\in Y_1\cap Y_2$. Conversely,
  if $y\in Y_1\cap Y_2$, then there exist $w_i\in \iM_i$ with $T_iw_i=y$
  for $i=1,2$.  Then $\bigl(\frac{w_1}{2},\frac{w_2}{2}\bigr)\in W$ and
  $T\bigl(\frac{w_1}{2},\frac{w_2}{2}\bigr)=y$.
\end{proof}

The following theorems refine some of the results
in~\cite{Arveson:67}, we will use them in the proof of
Theorem~\ref{t:inj-op-ran}.

\begin{theorem}\label{t:lin-ind}
  Let $\iA$ be $(n-1)$-transitive for some $n>1$ and
  $\iM=\overline{\iA\n\vec x}$ for some linearly independent $n$-tuple
  $\vec x=(x_1,\dots,x_n)$ in~$X^n$. Then either $\iM=X^n$ or $\iM$
  consists only of linearly independent $n$-tuples and zero.
\end{theorem}

\begin{proof}
  First, we show that if $\iN$ is a closed $\iA\n$-invariant subspace
  of $X^n$ such that $\iM\subseteq\iN$, then $\iN$ satisfies the
  following two properties.
  \begin{enumerate}
  \item\label{i:has-zero} If $\iN$ contains an $n$-tuple of the form
    $(u_1,\dots,u_{n-1},0)$ where $u_1,\dots,u_{n-1}$ are linearly
    independent, then $\iN=X^n$.
  \item\label{i:dep-indep} If $\iN$ contains both linearly independent
    and non-zero linearly dependent $n$-tuples, then $\iN=X^n$.
  \end{enumerate}
  Note that if $\iA$ is unital, then $\vec x\in\iM$; in this case
  applying~(\ref{i:dep-indep}) with $\iN=\iM$ would immediately yield
  the conclusion of the theorem.

  To prove~(\ref{i:has-zero}), notice that since $\iA$ is
  $(n-1)$-transitive then $X^{n-1}\oplus\{0\}\subseteq\iN$. Since
  $\iA$ is transitive, then $Ax_n\neq 0$ for some $A\in\iA$. It
  follows that there exists $(w_1,\dots,w_n)$ in $\iM$ with $w_n\neq
  0$. But $X^{n-1}\oplus\{0\}\subseteq\iN$ implies
  $(w_1,\dots,w_{n-1},0)\in\iN$, so that $(0,\dots,0,w_n)\in\iN$,
  which yields $\{0\}^{n-1}\oplus X\subseteq\iN$. Hence, $\iN=X^n$.

  To prove~(\ref{i:dep-indep}) suppose that $\iN$ contains a linearly
  independent $n$-tuple $\vec v=(v_1,\dots,v_n)$ and a non-zero
  linearly dependent $n$-tuple $\vec y=(y_1,\dots,y_n)$. 
  Without loss of generality,
  $y_1,\dots,y_k$ are linearly independent for some $k<n$, and
  $y_i=\sum_{j=1}^k\alpha_{ij}y_j$ as $i=k+1,\dots,n$. Since $\iA$ is
  $k$-transitive, there exists a sequence $(A_m)$ in $\iA$ such that
  $A_my_1\to v_1$ and $A_my_i\to 0$ as $i=2,\dots,k$. Then for $k+1\leq
  i\leq n$ we have $A_my_i\to\alpha_{i,1}v_1$. It follows that
  $$A_m\n\vec y\to
  (v_1,0,\dots,0,\alpha_{k+1,1}v_1,\dots\alpha_{n,1}v_1),$$
  so that the latter $n$-tuple belongs to~$\iN$. Subtracting it
  form $\vec v$ we see that $\iN$ contains an element of the form
  $(0,z_2,\dots,z_n)$, where $z_2,\dots,z_n$ are linearly independent.
  Now $\iN=X^n$ by~(\ref{i:has-zero}).
  
  To complete the proof, assume that $\iM$ contains a linearly
  dependent non-zero $n$-tuple and $\iM\neq X^n$. Notice that
  $\iM+[\vec x]$ is $\iA$-invariant, so that $\iM+[\vec x]=X^n$
  by~(\ref{i:dep-indep}).  In particular, there exists $\vec w\in\iM$
  and a scalar $\lambda$ such that $\vec w+\lambda\vec
  x=(x_1,0,\dots,0)$. Then $$\vec w=\bigl((1-\lambda)x_1,-\lambda
  x_2,\dots,-\lambda x_n)\bigr).$$
  If $\lambda=1$ then $\vec
  w=(0,-x_2,\dots,-x_n)$, so that $\iM=X^n$ by~(\ref{i:has-zero}),
  contradiction. If $\lambda\neq 1$ and $\lambda\neq 0$ then all the
  components of $\vec w$ are linearly independent, hence $\iM=X^n$
  by~(\ref{i:dep-indep}) which, again, would contradict our
  assumptions. It follows that $\lambda=0$, and so
  $(x_1,0,\dots,0)\in\iM$. Similarly,
  $(0,\dots,0,x_i,0,\dots,0)\in\iM$ for every $i=1,\dots,n$. It
  follows that $\vec x\in\iM$, so that $\iM=\iM+[\vec x]=X^n$,
  contradiction.
\end{proof}

\begin{remark}
  Assuming that $\iA$ is unital and $(n-1)$-transitive, it can be
  shown that any $\iM\in\Lat\iA\n$ is either in
  $\Lat\bigl(\iL(X)\bigr)\n$ (that is, consists of $n$-tuples
  satisfying a fixed set of linear dependence relations), or consists
  only of linearly independent $n$-tuples and zero. 
\end{remark}

\begin{corollary}\label{c:dep-dense}
  Let $\iA$ be transitive and assume that for any linearly independent
  $x_1$ and $x_2$ in $X$ there exists an operator $A\in\iA$ with
  $(Ax_1,Ax_2)$ non-zero and linearly dependent. Then $\iA$ is
  WOT-dense in $\iL(X)$.
\end{corollary}

\begin{proof}
  We will show by induction that $\iA$ is $n$-transitive for all $n\geq
  1$, then $\AWOT\!\!=\iL(X)$ by~\cite[Theorem~7.1]{Radjavi:73}. It is
  given that $\iA$ is 1-transitive. Suppose that $\iA$ is
  $(n-1)$-transitive for some $n>1$, and let $\vec x=(x_1,\dots,x_n)$ be a
  linearly independent $n$-tuple. Choose $A\in\iA$ so that
  $(Ax_1,Ax_2)$ is non-zero and linearly dependent. Then $\vec x$ is
  cyclic for $\iA\n$ by Theorem~\ref{t:lin-ind}. Hence, $\iA$ is
  $n$-transitive. 
\end{proof}

Our next result is a refinement of Corollary~2.5 in~\cite{Arveson:67}
(see also Lemma~8.8 in~\cite{Radjavi:73}).

\begin{theorem}[Graph Theorem]\label{t:graph}
  Suppose that $\iA$ is $n$-transitive. Then $\iA$ is not
  $(n+1)$-transitive if and only if there exists a closed operator
  $\vec T\colon\iD\subseteq X\to X^n$ where $\iD$ is a dense
  $\iA$-invariant linear subspace of $X$ and  $\vec
  T=T_1\oplus\dots\oplus T_n$, each $T_i$ commutes with $\iA$, and
  $(x,T_1x,\dots,T_nx)$ is linearly independent for each non-zero
  $x\in\iD$. In particular, each $T_i$ is one-to-one and
  non-scalar.
\end{theorem}

\begin{proof}
  Suppose that $\iA$ is $n$-transitive. If $\iA$ is not
  $(n+1)$-transitive then by Theorem~\ref{t:lin-ind} there exists a
  closed $\iA^{(n+1)}$-invariant subspace $\iM$ of $X^{n+1}$ such that
  every non-zero element of $\iM$ is a linearly independent
  $(n+1)$-tuple.  We claim that $\iM$ is the graph of a closed
  operator $\vec T\colon\iD\subseteq X\to X^n$ satisfying the required
  conditions.  Indeed, let
  $\iD=\bigl\{x_0\mid(x_0,x_1,\dots,x_n)\in\iM\bigr\}$. For
  $(x_0,x_1,\dots,x_n)\in\iM$ define $T_ix_0=x_i$ for $i=1,\dots,n$,
  and put $\vec T=T_1\oplus\dots\oplus T_n$. Notice that $\vec T$ is
  well-defined: suppose that $(x_0,x_1,\dots,x_n)$ and
  $(x_0,x'_1,\dots,x'_n)$ are both in~$\iM$, then
  $(0,x_1-x'_1,\dots,x_n-x'_n)\in\iM$, but this $(n+1)$-tuple is
  linearly dependent, hence equals zero, so that $x_i=x'_i$ for all
  $1\leq i\leq n$.  Clearly, $\vec T$ is closed because $\iM$ is closed.
  
  Now let $x\in\iD$, then $(x,T_1x,\dots,T_nx)\in\iM$. Since $\iM$ is
  $\iA^{(n+1)}$-in\-var\-iant, if $A\in\iA$ then
  $(Ax,AT_1x,\dots,AT_nx)\in\iM$. It follows that $Ax\in\iD$, so that
  $\iD$ is $\iA$-invariant and, therefore, dense in~$X$. Furthermore,
  it also follows that $AT_ix=T_i(Ax)$ for all $i=1,\dots,n$, so that
  $T_i$ commutes with~$\iA$. Finally, $(x,T_1x,\dots,T_nx)$ is in
  $\iM$ for each non-zero $x\in\iD$, so that this $(n+1)$-tuple is
  linearly independent. It follows that each $T_i$ is non-scalar.  It
  also follows that $T_ix=0$ implies $x=0$, so that $T_i$ is
  one-to-one.
  
  To prove the converse, suppose that $\vec T$ is such an operator,
  and let $\iM$ be the graph of $\vec T$. Then $\iM$ is closed and
  every non-zero element of $\iM$ is linearly independent. It follows
  that $\iM\neq X^{n+1}$. It is easy to see that $\iM$ is
  $\iA^{(n+1)}$-invariant, so that every linearly independent
  $(n+1)$-tuple in $\iM$ is not cyclic for $\iA^{(n+1)}\!\!$, hence $\iA$
  is not $(n+1)$-transitive.
\end{proof}

Now we are ready to prove Theorem~\ref{t:inj-op-ran}.

\begin{proof}[Proof of~Theorem~\ref{t:inj-op-ran}]
  Suppose that $\iA$ has no non-trivial invariant injective operator
  ranges of order $n$ for some $n$. It follows that $\iA$ has no
  non-trivial invariant injective operator ranges of order $k$
  whenever $1\leq k\leq n$. In particular, $\iA$ is transitive.
  
  Note first, that $\iA'=[I]$.  Indeed, let $T\in\iA'$. Pick
  $\lambda\in\sigma(T)$, and show that $T=\lambda I$. Suppose not,
  then $\Range(T-\lambda I)\neq 0$. Since $\Null(T-\lambda I)$ is an
  $\iA$-invariant closed subspace then $T-\lambda I$ is one-to-one.
  Then $\Range(T-\lambda I)$ is an $\iA$-invariant injective operator
  range of order 1, so that $\Range(T-\lambda I)=X$. But this would
  mean that $T-\lambda I$ is invertible, contradiction.
  
  Suppose that $\iA$ is not $n$-transitive, and show that this leads
  to a contradiction.  Let $k$ be the greatest such that $\iA$ is
  $k$-transitive but not $(k+1)$-transitive, then $1\leq k<n$.  Let
  $\vec T\colon\iD\to X^k$ where $T=T_1\oplus\dots\oplus T_k$ is as in
  Theorem~\ref{t:graph}. Let $\iM$ be the graph of~$\vec T$.  Also,
  let $P\colon\iM\to X$ be the projection on the first component. Then
  $\iD=\Range P$ is an injective operator range of order $(k+1)$. It
  follows that $\iD=X$, so that $T_i\in\iL(X)$ for each $i=1,\ldots,k$
  and, therefore, $T_i\in\iA'$.  But then $T_i$ have to be scalar
  because $\iA'=[I]$, contradiction.
\end{proof}


\begin{thebibliography}{DHP01}

\bibitem[Arv67]{Arveson:67}
William~B. Arveson.
\newblock A density theorem for operator algebras.
\newblock {\em Duke Math. J.}, 34:635--647, 1967.

\bibitem[DHP01]{Donsig:01}
Allan Donsig, Alan Hopenwasser, and David~R. Pitts.
\newblock Automatic closure of invariant linear manifolds for operator
  algebras.
\newblock {\em Illinois J. Math.}, 45(3):787--802, 2001.

\bibitem[Foi72]{Foias:72}
Ciprian Foia{\c{s}}.
\newblock Invariant para-closed subspaces.
\newblock {\em Indiana Univ. Math. J.}, 21:887--906, 1971/72.

\bibitem[Hop01]{Hopenwasser:01}
Alan Hopenwasser.
\newblock Invariant linear manifolds for {CSL}-algebras and nest algebras.
\newblock {\em Proc. Amer. Math. Soc.}, 129(2):389--395, 2001.

\bibitem[Lam71]{Lambert:71pacific}
Alan Lambert.
\newblock Strictly cyclic operator algebras.
\newblock {\em Pacific J. Math.}, 39:717--726, 1971.

\bibitem[Ric50]{Rickart:50}
C.~E. Rickart.
\newblock The uniqueness of norm problem in {B}anach algebras.
\newblock {\em Ann. of Math. (2)}, 51:615--628, 1950.

\bibitem[RR73]{Radjavi:73}
Heydar Radjavi and Peter Rosenthal.
\newblock {\em Invariant subspaces}.
\newblock Springer-Verlag, New York, 1973.
\newblock Ergebnisse der Mathematik und ihrer Grenzgebiete, Band 77.

\bibitem[Yood49]{Yood:49}
Bertram Yood.
\newblock Additive groups and linear manifolds of transformations between
  {B}anach spaces.
\newblock {\em Amer. J. Math.}, 71:663--677, 1949.

\end{thebibliography}
\end{document}